\documentclass{lms}
\usepackage{amsmath,amsfonts,amscd,amssymb,latexsym,comment}
\usepackage{enumerate,graphicx,psfrag,subfigure,changebar,oldgerm}
\title{Centraliser Dimension and Universal Classes of Groups}
\author{A.J. Duncan, 
 I.V. Kazatchkov
 \and V.N. Remeslennikov}
\classno{20E15 (primary), 20F10, 03B25 (secondary)}
\extraline{This research was carried out with the support of EPSRC
  grant GR/S71200, ``Algebraic geometry over groups and algorithms for
  Artin groups''.}
\def\nul{\emptyset }

\def\b{\beta }

\def\i{\iota }

\def\G{\Gamma }
\def\g{\gamma }
\def\a{\alpha }

\def\fC{{\textswab C}}
\newtheorem{thm}{Theorem}[section]

\newtheorem{cor}[thm]{Corollary}
\newtheorem{prop}[thm]{Proposition}
\newnumbered{ax}[thm]{Axiom}
\newnumbered{defn}[thm]{Definition}
\newnumbered{expl}[thm]{Example}
\newnumbered{rem}[thm]{Remark}
\newnumbered{con}[thm]{Conjecture}
\newnumbered{que}[thm]{Question}
\newnumbered{prob}[thm]{Problem}
\numberwithin{equation}{section}
\newcommand{\CT}{\texttt{CT}}
\newcommand{\cd}{\texttt{cdim}}
\newcommand{\DC}{\texttt{CD}}
\newcommand{\ttz}{\texttt{z}}
\newcommand{\CSA}{\texttt{CSA}}
\newcommand{\US}{\texttt{US}}
\newcommand{\Cyc}{\texttt{Cyc}}
\newcommand{\Int}{\texttt{Path}}

\newcommand{\dis}{\texttt{Dis}}
\newcommand{\ldis}{\texttt{LDis}}
\newcommand{\ZZ}{\ensuremath{\mathbb{Z}}}

\newcommand{\NN}{\ensuremath{\mathbb{N}}}

\newcommand{\cK}{\mathcal{K}}

\newcommand{\Th}{\textsf{Th}}

\newcommand{\ucl}{\texttt{ucl}}

\newcommand{\la}{\langle}
\newcommand{\ra}{\rangle}
\newcommand{\sdc}{>}
\newcommand{\sac}{<}

\newcommand{\edc}{\ge}
\newcommand{\eac}{\le}

\newcommand{\mbf}{\mathbf}
\newcommand{\be}{\begin{enumerate}}
\newcommand{\ee}{\end{enumerate}}
\begin{document}
\maketitle
\begin{abstract}
In this paper we establish results that will be required
for the study of the algebraic geometry of partially-commutative
groups. We define classes of groups axiomatised by sentences
determined by a graph. Among the classes which arise this way we
find $\CSA$ and $\CT$ groups. We study the centraliser dimension
of a group, with  particular attention to the height of the lattice
of centralisers, which we call the centraliser dimension of the group.
The behaviour of centraliser dimension under several common group
operations
is described. Groups with centraliser dimension $2$ are studied in
detail. It is shown that $\CT$-groups are precisely 
those with centraliser dimension $2$ and  trivial centre.
\end{abstract}
\section{Introduction}
The purpose of this paper is to lay foundations for the study of 
equations over groups and in particular over free
partially-commutative
groups. We construct universal and existential sentences based on
graphs
and relate these to groups. 
The formula $\phi(\Gamma)$ which we introduce, given a graph $\Gamma$,
and the properties developed below 
suggest the following general question.
\begin{que}
Let $\Gamma_1$ and $\Gamma_2$ be two finite connected graphs and
suppose that the formula $\phi(\Gamma_1)$ is logically equivalent
to the formula $\phi(\Gamma_2)$ for all groups from some class $\mathcal K$.
What can be said about $\Gamma_1$ and $\Gamma_2$?
\end{que}
We also investigate properties of centralisers of groups.
Among other things we show that the class of groups which has a centraliser
lattice of finite height $m$ is universally axiomatisable and
describe the behaviour of this class under various group operations.
In subsequent papers we plan to apply the results of this
work to find the centraliser dimension of free
partially-commutative groups and to investigate the problem of
universal equivalence for this class of groups. 
Our interest in
these problems is inspired by the importance of such results
in 
algebraic geometry over groups  (see \cite{alg1}, \cite{alg2}): and in 
particular over 
free partially-commutative
groups.

We begin by considering classes of
groups axiomatised by certain sentences in the first order
language corresponding to graphs. In this way, we arrive at
certain classes of groups, some new, among which are the
well-known classes of $\CT$- and $\CSA$-groups. 

In Section \ref{section:cd} we turn to the study  of what we
call the centraliser dimension of a group. This coincides with 
the notion of height of the centraliser lattice of a group, introduced
by R. Schmidt \cite{Schmidt94}. The lattice of centralisers of  various
groups, have been investigated by numerous authors: see for example
\cite{Ito},  
\cite{Zaleski65}, \cite{Vas}, \cite{LR}, \cite{Schmidt95}, \cite{Weh}, 
\cite{Kos73}, \cite{kegel73}, \cite{BryantHartley79}, \cite{Bryant79},
\cite{AS}, \cite{Bl}, \cite{Wag} and
\cite{MS}.
In particular, a detailed account of results in the field
can be found in V.~A.~Antonov's book  \cite{book}.
Here we show that the groups which have centraliser lattice of finite 
height are universally axiomatisable.

Next we investigate the behavior of the centraliser dimension under
several group operations: namely free products, direct products
and free products with amalgamation by their centres. We also 
study of groups with centraliser
dimension $2$. The groups with trivial centre which are of centraliser
dimension $2$ are shown to coincide with the class of $\CT$-groups.
Examples show that when the centre of the group is non-trivial the 
picture is far more complex.
\section{Universal classes and some notions from model theory}
\label{section:uc}
\subsection{Preliminaries}

We recall here some basic notions of model theory that we require.
For more details we refer the reader to \cite{CK}.
The standard language of group theory, which we denote by $L$,
consists of a symbol for multiplication `$\cdot$', a symbol for
inversion $^{-1}$, and a constant symbol for the identity. We take
$X=\left\{ x_1,x_2,  \ldots \right\}$ as the set of variables of our
language and define
$X^{-1}=\left\{x^{-1}|x \in X \right\}$ and $X^{\pm 1}=X \cup
X^{-1}$. A {\em term} is an element of the free semigroup on
$X^{\pm 1}$. An {\em atomic formula} is an expression of the form
$w=1$, where $w$ is  a term. A {\em formula} in $L$ is either an
atomic formula or one of $\theta \lor \phi$, $\theta \land \phi$,
$\neg \phi$, $\forall x \phi$ or $\exists x \phi$, where $\theta$ and
$\phi$ are formulas (and $\lor$, $\land$, $\neg$, $\forall$ and
$\exists$ have their usual meanings). If $\theta$ is a formula
and $S$ is a subset of $X$ then we write $\theta(S)$ to
indicate that the variables which occur in $\theta$ are all elements
of
$S$.
It follows from standard first order logic that any formula is
logically
equivalent to a formula of the type
\begin{equation} \label{eq:logf}
Q_1y_1Q_2y_1\cdots Q_my_m \psi(x_1, \ldots ,x_m,y_1,\ldots ,y_n),
\end{equation}
where $Q_i \in \left\{ \forall, \exists \right\}$
and $\psi(x_1, \ldots ,x_m,y_1,\ldots ,y_n)$ is a formula.
We shall therefore assume formulas have this form.
Those of
the $x_i$'s occuring  are called {\em free} variables.
If \eqref{eq:logf} has no free variables it is
called
a {\em sentence} in $L$.

Let $G$ be a group. We assume that, for a sentence $\phi$ in $L$, the meaning
of ``$\phi$ holds in $G$'' is understood. For example the sentence
$\forall x \forall y ([x,y]=1)$ holds in $G$ if and only if $G$ is
Abelian.
If $\phi(x_1,\ldots,x_m)$ is an arbitrary formula in $L$ then
we denote by $\phi(g_1,\ldots ,g_m)$ the element of $G$ obtained by
substiting $g_i$ for $x_i$ in $\phi$, whenever $x_i$ is a free
variable of $\phi$. Let $g$ be the sequence $g=g_1,\ldots ,g_m\in G^m$.
We write $G\models \phi(g)$ if $\phi(g)$ holds in $G$.
For example, let $\phi(x_1,x_2,x_3)$ be the formula
$\forall x_1\exists x_2 (x_1\cdot x_2=x_3)$ and let $g=g_1,g_2,g_3\in G^3$.
Then $\phi(g)$ is $\forall x_1\exists x_2 (x_1\cdot x_2=g_3)$ so
$G\models \phi(g)$.
The {\em truth domain} of $\phi$ over $G$ is
\[
\phi(G)= \left\{g \in G^m| G \models \phi(g) \right\}.
\]
If $\phi(G)=G^m$ then we write $G\models \phi$ and say that $\phi$ is
{\em satisfied} by $G$, $\phi$ is {\em valid} in $G$, $\phi$ {\em holds}
in $G$ or that $G$ is a $\phi${\em -group}.
Of course, since a sentence has no free variables, when $\phi$ is a
sentence,
this reduces to the
notion of ``holds in $G$'' that we assumed above.
Let $\cK$ be a class of groups. Then we say that $\theta$ and
$\phi$ are {\em logically equivalent} in $\cK$ if $G\models \theta$ if
and only if $G\models \phi$, for all groups $G$ from $\cK$.  We say
that
$\cK$ is {\em axiomatisable} by a set of sentences $S$ if
$\cK$ consists of all groups $G$ such that $G\models s$, for all $s\in
S$.

If $G$ is a group then the set $\Th_\forall (G)$ of all
\emph{universal sentences} (i.e. $Q_1=\ldots = Q_m= \forall$ in
Formula (\ref{eq:logf})) which are valid in $G$ is called the
\emph{universal theory} of $G$. By the definition, two groups $G$
and $H$ are \emph{universally equivalent} if $\Th_\forall
(G)=\Th_\forall (H)$, in which case we write $G\equiv_\forall H$.
The \emph{universal closure} $\ucl (G)$ of a group $G$ consists of
all groups $H$ such that $\Th_\forall (G)\subseteq \Th_\forall
(H)$. A class of groups $\mathcal K$ is {\em universally
axiomatisable} if it can be axiomatised by a set of universal
sentences.
 The {\em existential theory} $Th_\exists(G)$ of $G$ is defined
analagously, as are {\em existential equivalence} and {\em existential closure}.
Notice that conditions $G\equiv_\forall H$ and $G \equiv
_\exists H$ are equivalent.

Let $G$ be a group and $M$ be a set of elements of $G$. Then the
set $M$ together with induced partial group operation on it is called
a {\em partial model} of $G$. On
the set of partial models of $G$ the notion of isomorphism of
partial models arises naturally. 
The
following Proposition follows from well-known facts of model theory.

\begin{prop} \label{prop:pmodel}
Let $G$ and  $H$ be groups. Then $G\equiv_\forall H$ if and only if  every
finite partial model of $G$ is isomorphic to a finite
partial model of $H$, and vice-versa.
\end{prop}

\begin{cor}
Free non-Abelian groups (of arbitrary rank) are universally
equivalent.
\end{cor}

\begin{defn} \label{defn:I-discr}
Let $G$ and $H$ be groups. We say that $G$ is
{\em discriminated} by $H$ if, for every finite subset $\left\{
g_1 ,\ldots,g_m \right\}$ of non-trivial elements of $G$,
there exists a homomorphism $\varphi: G \rightarrow H$ such that
$\varphi (g_i ) \ne 1$, for $i = 1, \ldots, m$. The set of
all groups discriminated by $H$ is denoted $\dis(H)$.
\end{defn}

If every finitely generated subgroup of $G$ is discriminated by
$H$ then we say that $G$ is {\em locally discriminated} by $H$.
The set of all groups locally discriminated by $H$ is denoted
$\ldis(H)$. The next Proposition follows immediately from
Proposition \ref{prop:pmodel} and the fact that if $G$ is locally
discriminated by $H$ then, for every finite subset $\left\{
g_1 ,\ldots,g_m \right\}$ of non-trivial elements of $G$, we may choose
$\varphi: G \rightarrow H$ such that $\varphi(g_i)\neq \varphi(g_j)$,
when $i\neq j$.

\begin{prop} \label{prop:ldis}
Suppose that $G$ is locally discriminated by $H$ and that $H$ is
locally discriminated by $G$. Then $G \equiv_\forall H$.
\end{prop}

\begin{cor}
Arbitrary non-trivial torsion-free Abelian groups are
universally equivalent.
\end{cor}
\begin{proof} An infinite cyclic group is discriminated by any
torsion-free Abelian group. Hence, by Proposition \ref{prop:ldis}, it suffices to prove that
every Abelian group of finite rank is discriminated by an infinite
cyclic group, and this is easily verified.  \end{proof}

\subsection{Logical Formulas and Universal Classes}\label{s:logic}

We next describe some important logical group formulas in the language $L$, involving
 commutation of group elements. Our conventions are that
$\left[ x,y\right]=x^{-1}y^{-1}xy$ and $x^y=y^{-1}xy$.

Let $\g(x_1, \ldots, x_k,y)$ denote the formula
\[
\bigwedge \limits_{i=1}^k \left[x_i,y\right]=1.
\]
A $(k+1)$-tuple $(g_1, \ldots, g_k, g) \in  G^{k+1}$ is in the
truth domain $\g(G)$ of $\g$  if and only if $g$ is contained in
the centraliser $C_G(g_1, \ldots, g_k)$ of $\left\{ g_1,\ldots, g_k
\right\}$ in $G$ (see Section \ref{s:defn}).
Thus it is natural to use  $y\in C(x_1,\ldots, x_k)$ to denote
$\g(x_1, \ldots, x_k,y)$.
Similarly we use $y \notin C(x_1,\ldots, x_k)$ for the negation of
$\g$.
Similarly, by $x \in Z$ we denote the formula
$\forall y\left[y, x\right]=1$, since the truth domain of this
formula over a group $G$ is the centre $Z(G)$ of $G$.

The {\em commutativity axiom} is the sentence
$\forall x,y \left[ x,y \right]=1$: valid in  the group $G$ if and only if  $G$
is Abelian.
The {\em commutative transitivity or $\CT$ axiom} is the sentence $\CT(x,y,z)$
given by
\[
\forall x, y, z \ (x \ne 1 \land x\in C(y,z) \rightarrow \left[ y, z \right]=1).
\]
%If this axiom is satisfied by $G$ then we say that $G$ is a
%$\CT$-group.
The $\CT$ axiom  is logically equivalent (in the class of all groups)
to the sentence
\begin{gather*}
\forall x,y,z(x=1\lor y=1\lor z=1\lor x=y\lor x=z \lor y=z\\\lor
[x,y]\neq 1\lor [x,z]\neq 1 \lor [y,z]=1).
\end{gather*}
Thus $G$ is
$\CT$-group if and only if the centraliser
of every nontrivial element of $G$ is an Abelian subgroup; which must therefore be
a maximal Abelian subgroup.

The $\CSA$ {\em axiom} is the sentence $\CSA(x,y,t)$ given by
\[
\forall x,y,t(x\neq 1 \land x\in C(y)\land x^t\in C(y) \rightarrow t
\in C(y)).
\]
%If $G$ satisfies the $\CSA$ axiom then it is said to be a $\CSA$-group.
%%%%%%%%%
A subgroup
$M$ of a group $G$ is {\em conjugacy-separable} or {\em malnormal} 
if $M\cap M^g=1$, for all $g\in G\backslash M$. Thus a group is a $\CSA$-group
if and only if all centralisers of single elements are conjugacy-separable.

Clearly a $\CSA$-group is a $\CT$-group. However the converse does not hold.
The free product of two cyclic
groups of order two is an example of a group which is a
$\CT$-group but not a $\CSA$-group. In fact, if the factors
are generated by $a$ and $b$ then the centraliser of $ab$ does not
contain $b$  but is fixed under conjugation by $b$ and hence is not
conjugacy-separable. Therefore, the
class of $\CT$-groups is wider than the class of
$\CSA$-groups.

A $\CSA$-group also satisfies the following which we call the 
{\em unilateral-separability} or $\US$-axiom. The $\US$-axiom is the 
sentence $\US (x,y)$ given by
\[\forall x, y \ (x^y\in C(x) \rightarrow y\in C(x)).\]
%As usual, groups satisfying this axiom are called $\US$-groups.
Again the class of $\US$-groups is wider than the class of $\CSA$-groups. 
For example if $F$ is a free group of rank
$2$ and $C$ is infinite cyclic then it is easy to see that $F\times C$ is 
a $\US$-group. However $F\times C$ is not a $\CT$-group so is not a $\CSA$-group.

Now suppose that $G$ is both a $\CT$-group and a $\US$-group. Let $x,y$ and $z$ 
be elements of $G$ with $x\neq 1$ and both $x\in C(y)$ and $x^z\in C(y)$. Then,
as $G$ is a $\CT$-group, $x^z\in C(x)$ and, as $G$ is a $\US$-group, $z\in C(x)$.
Using the $\CT$-axiom again $z\in C(y)$, so $G$ is a $\CSA$-group. Hence the 
$\CSA$-axiom is logically equivalent, in the class of all groups, to the
sentence
\[
\forall x, y, z\ (\CT(x,y,z) \land \US(x,y)).
 \]
%%%%%%%%%%
\begin{comment}
Clearly a $\CSA$-group is a $\CT$-group.
Hence the $\CSA$ axiom is logically equivalent, in the class of all groups, to the
sentence
\[
\forall x, y, z, u,t \ (\CT(x,y,z) \land x \in C(u) \land
x^t \in C(u) \rightarrow t \in C(u)).
\]
It follows that $G$ is a $\CSA$-group if and only if it is a
$\CT$-group in which the centralisers of non-trivial elements are
malnormal. (A subgroup
$M$ of a group $G$ is malnormal if $M\cap M^g=1$, for all $g\in G\backslash M$.)
\end{comment}

It is also not hard to show that a group is $\CSA$ if and only if all
maximal Abelian subgroups are conjugacy-separable.
For more details of $\CSA$ groups see \cite{MR2}.

\subsection{Paths and Cycles}

We describe a number of existential sentences which
encapsulate the relations of commutativity of a finite set of
elements. These sentences for commutativity relations
are indexed by certain graphs and so we shall begin by defining
formulas $\theta(\G)$ and  $\phi (\Gamma)$
%%%%%%%%%%%%%%%%%%%%%%%%%
corresponding to an arbitrary graph
$\Gamma$. Let $\Gamma$ be a  graph with vertices $V(\Gamma)$ and
edges $E(\Gamma)$. For notational simplicity we
assume that $V(\G)$ is a subset of $X$ and that $V(\G)=\{x_1,\ldots,x_n\}$.
The sentence $\phi(\Gamma)$ is defined to be
$\exists x_1 \ldots \exists x_n \theta(\G)$, where
$\theta(\G)$ is
the conjunction
of the following formulas.
\begin{enumerate}
\item
$x \ne 1$, for all $x \in V(\Gamma)$;
\item
$x \ne y$, if $x$ and $y$ are disjoint vertices of $\G$;
\item
$[x,y]=1$, whenever there is an edge in
$\Gamma$ connecting $x$ and $y$ and
\item $[x,y]\ne 1$ if there
is no such edge.
\ee
Thus $G\models \phi(\G)$ if and only if $G\models
\theta(\G)(g_1,\ldots, g_n)$, for some $g_1,\ldots ,g_n\in G^n$. In
this case we call the sequence $g_1,\ldots ,g_n$ an {\em implementation } of $\G$ in $G$ and
say that $G$ {\em admits} the graph $\G$.
%~\\
Let $\Phi(\G)$ be the class of all groups
   in  which the sentence $\phi (\Gamma)$ is satisfied and let $\Phi(\neg
    \Gamma)$ be the complement of the class
    $\Phi(\G)$. Clearly,  since $\neg
    \phi(\Gamma)$ is a universal formula, the class $\Phi(\neg
    \G)$ is a universal class.
%%%%%%%%%%%%%%%%%%%%%%%5

The {\em path graph} $ \Int_l$ of length $l$ is a tree with $l+1$ vertices precisely
two of which have degree one. Our first family of sentences for commutativity relations
is indexed by the
path graphs of positive length.%\\[1em]

The {\em length-one-path axiom} is the sentence $\phi(\Int_1)$, that is
\[
\exists x_1, x_2 \ (x_1 \ne 1 \land x_2\neq 1 \land x_1 \neq x_2 \land \left[ x_1, x_2
\right]=1).
\]
The negation of this sentence $\neg \phi(\Int_1)$
is
\[
\forall x_1, x_2 \ (x_1 = 1 \lor x_2 =1 \lor x_1=x_2 \lor \left[ x_1, x_2
\right]\neq 1).
\]
Clearly this sentence is satisfied by groups of order at most $2$. If
$G$ is a group of order more than $2$ then either $G$ has an  element $g$ of
order $3$ or more, or all non-trivial elements of $G$ have order $2$. In the
former case $\neg \phi(\Int_1)$ does not hold in $G$ since we may take
$x_1=g$ and $x_2=g^2$.
In the latter case, since the order of $G$ is more than $2$, it has non-trivial elements
$a$ and $b$ with $a\neq b$ and $[a,b]=1$, so $\neg \phi(\Int_1)$ does not hold.
Therefore $\Phi(\neg \Int_1)$ consists of
of the trivial group and the cyclic group of order $2$.
It follows that $\neg \phi(\Int_1)$
is logically equivalent, in the class
of all groups,  to the universal sentence
\[
\forall x_1, x_2 \ (x_1 = 1 \lor x_2 =1 \lor x_1=x_2).
\]

  %  \\
The {\em length-two-path axiom} is the sentence $\phi(\Int_2)$ given by
\begin{gather} \notag
  \begin{split}
    \exists x_1, x_2,x_3 & \ (\bigwedge\limits_{i=1}^{3} x_i \ne 1
    \bigwedge_{\substack{ i,j =1\\i\ne j}}^{3}(x_i \ne x_j) \wedge \\
    \wedge & \left[x_1, x_2\right]=1 \wedge [x_2,x_3]=1 \wedge
    [x_1,x_3] \ne 1).
  \end{split}
\end{gather}
(Assuming that $x_2$ is the vertex of $\Int_2$ of degree $2$.)
Negation of  $\phi(\Int_2)$ is a universal sentence,
$\neg \phi(\Int_2)$,
logically equivalent, in the class of all groups,
to the $\CT$ axiom.
Therefore $\Phi(\neg\Int_2)$ is
is the class of
$\CT$-groups.

Similarly, the  {\em length-}$l${\em -path axiom}, for $l \ge 3$, is defined to
be $\phi(\Int_l)$.
%In the next section we shall give an example of a group which
%satisfies  $\phi(\Int_2)$ but does not satisfy $\phi (\Int_3)$.\\
Now $\Phi(\Int_l)$,  the class of groups which satisfy
$\phi(\Int_l)$, clearly satisfies
$\Phi(\Int_l)\ge\Phi(\Int_{l+1})$, for all $l\ge 1$.
It is easy to see that if $F$ is a free group of rank $2$ and $C$ is
infinite cyclic then $F\times C$ satisfies  $\phi(\Int_2)$
 but does not satisfy $\phi (\Int_3)$. 
Moreover Blatherwick \cite{B} has shown that 
in fact $\Phi(\Int_l)>\Phi(\Int_{l+1})$, for $l\ge 1$. Thus we have the 
following chain of inclusions.
\[
\Phi(\Int_1)\sdc \Phi(\Int_2)\sdc \Phi(\Int_3) \sdc \Phi(\Int_4) \sdc
\Phi(\Int_5) \sdc \cdots .
\]

We next consider a family  of existential sentences indexed by
cycle graphs. The {\em cycle graph}, or $l${\em -cycle},
$\Cyc_l$ is the connected graph with $l$
vertices which is regular of degree $2$.
The {\em three-cycle axiom} is the sentence $\phi(\Cyc_3)$
\begin{multline*}
  \exists x_1, x_2,x_3 \
  (
    \bigwedge\limits_{i=1}^{3} x_i \ne 1
    \bigwedge_{\substack{i,j =1\\i\ne j}}^{3}(x_i \ne x_j)
    \wedge \\
    \wedge \left[x_1, x_2\right]=1 \wedge \left[x_2,x_3\right]=1 \wedge
    \left[x_1,x_3\right] = 1
  ).
\end{multline*}
\begin{prop}\label{prop:cyc3}
The class $\Phi(\Cyc_3)$ consists of all groups except those
of order less than $4$ and the dihedral group $D_6$.
\end{prop}
\begin{proof}
Clearly $\phi(\Cyc_3)$ does not hold in any group of order less than $4$ or
in $D_6$. Suppose that $G$ is a group of order at least $4$ and that $G\ncong D_6$.
We shall find a sequence $\mbf v=a,b,c$ of elements of $G$ which
satisfy
$\theta(\Cyc_3)(x_1,x_2,x_3)$.
If there are non-trivial elements $a$ and $b$ of $G$ such that $a\neq
b^{\pm 1}$ and $[a,b]=1$ then we may take $\mbf v=a,b,ab$. Hence we
may assume that $G$ contains no such elements.
If $G$ has an element $a$ of order $4$ or more then
we can take $\mbf v= a, a^2, a^3$. Assume then that all elements of $G$ are of
order at most $3$. If $G$ has no element of order $3$ then $G$ has
distinct non-trivial elements $a$ and $b$ with $a\neq b^{\pm 1}$ and
$[a,b]=1$, a contradiction.
Thus we may assume $G$ has elements of order $3$. Suppose that $G$ has
no elements of order $2$. Then we may choose non-trivial elements
$a$ and $b$ of $G$
with $a\neq b^{\pm 1}$ and $[a,b]\neq 1$.
Since $a$, $b$, $ab$, $ab^2$ and $a^2b^2$ have order $3$,
it follows that
$[a,a^b]=[a^2, a^b]=1$. It is then easy to check that
we may take $\mbf v=a, a^2,a^b$.
This leaves the case where $G$ has elements
of order $2$ and $3$. Let $a$ and $b$ be elements of $G$ with $|a|=2$,
$|b|=3$ and $[a,b]\neq 1$. Suppose first that $|ab|=3$. Then
$aba=b^2ab^2$ and so $[a,a^b]=ab^2abab^2ab=(ab)^3=1$.
Similarly $[a,a^{b^2}]=1$ and a straightforward check shows that
$\mbf v=a, a^b,a^{b^2}$ satisfies $\theta(\Cyc_3)$. Now suppose that $|ab|=2$ (with
$|a|=2$, $|b|=3$ and  $[a,b]\neq 1$ as before). Then $\la a,b\ra\cong D_6$. As
$G\ncong D_6$ it must contain an element $c\notin \la a,b\ra$.
Then, given our intial assumptions, we have
$[a,c]\neq 1$ and $[b,c]\neq 1$. If $|c|=3$ and $\la a,c\ra\ncong D_6$ then
we can find $\mbf v$ satisfying $\theta(\Cyc_3)$,  as above. If $\la a,c\ra \cong D_6$ then
$|ac|=2$. As $ac\notin \la a, b\ra$ we may replace $c$ with $ac$. Thus, without
loss of generality, we may assume that $|c|=2$ and that $\la a,b\ra\cong \la c,b\ra
\cong D_6$. If $|ac|=2$ then we have $a\neq c^{\pm 1}$ and $[a,c]=1$,
a contradiction. Hence
$|ac|=3$. Then
$(acb)^3=acbacbacb=ab^2cacab^2cb=b(ac)^3b^2=1$. Similarly
$(acb^2)^3=((ac)^2b^2)^3=1$ and, as in the case where $G$ has no
element of order $2$, we may now take
$\mbf v=ac, (ac)^2, (ac)^b$.
\end{proof}

The {\em four-cycle axiom} is the sentence $\phi(\Cyc_4)$ given by
\begin{multline*}
    \exists x_1, x_2,x_3,x_4 \ (\bigwedge\limits_{i=1}^{4}
    x_i \ne 1   \bigwedge_{\substack{i,j=1\\i\ne j}}^{4}(x_i \ne x_j)
    \wedge \left[x_1, x_2\right]=1 \wedge\\
    \left[x_1, x_4\right]=1 \wedge [x_2,x_3]=1 \wedge
    [x_3,x_4] = 1\wedge [x_1,x_3]\ne 1  \wedge
    [x_2,x_4] \ne 1).
\end{multline*}

If we identify $x_1,x_2,x_3,x_4$ with an implementation of $\Cyc_4$ in
$G$ then in the drawing of
$\Cyc_4$ below letters
connected by an edge commute and letters that are not connected by
an edge do not.
\begin{center}
\psfrag{a}{$x_1$}
\psfrag{b}{$x_2$}
\psfrag{c}{$x_3$}
\psfrag{d}{$x_4$}
\includegraphics[scale=0.3]{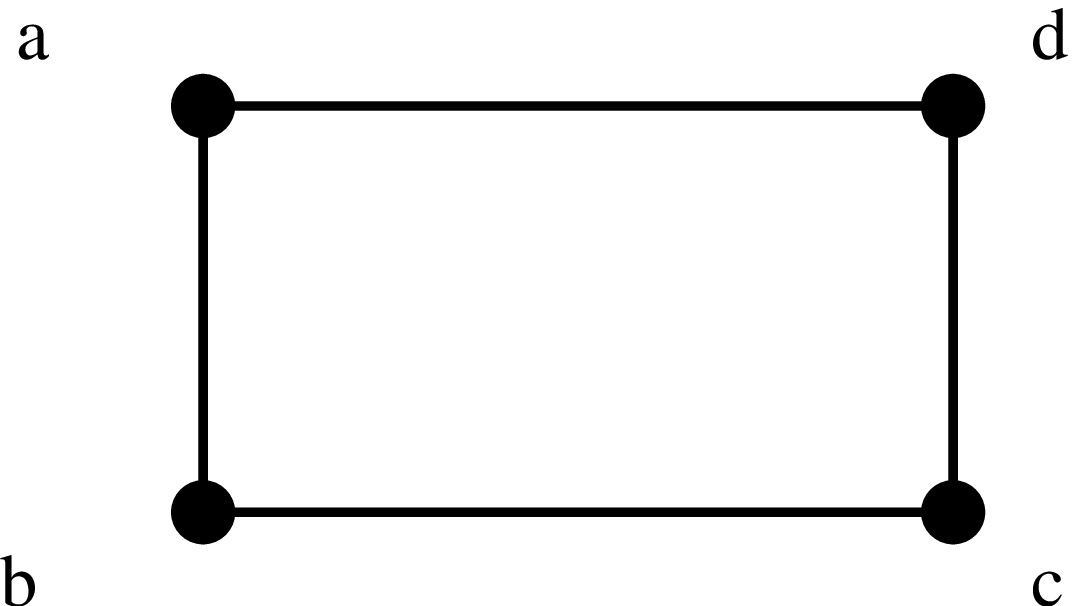}
\end{center}
The negation of $\phi(\Cyc_4)$ is a universal
sentence which is satisfied by $G$ if and only if $G$ admits no
$4$-cycle.  It is clear that if $G$ is a $\CT$-group then it
admits no $4$-cycle. 
However the converse does not hold: for example $D_8$, the dihedral group
of order $8$, admits no $4$-cycle
and is
not a
$\CT$-group. This means that $D_8*D_8$ admits no $4$-cycle, is not
a $\CT$-group and in addition has trivial centre.

Similarly the $5${\em -cycle axiom} is the sentence $\phi(\Cyc_l)$, for $l \ge 5$.
Blatherwick \cite{B} has shown that, although $\Phi(\Cyc_3)>\Phi(\Cyc_4)$, for
$n\ge 5$ and $m=n \pm 1$, $\Phi(\Cyc_n)\ngeq \Phi(\Cyc_m)$.
%

\begin{comment}
\begin{thm}
If $\Phi(\Gamma_1)$ is logically equivalent to $\Phi (\Gamma_2)$
in the class of nilpotent groups of class $2$, where $\Gamma_1$ and $\Gamma_2$
are connected
graphs, then the graph  $\Gamma_1$  is
logically equivalent to the graph $\Gamma_2$.
\end{thm}
\end{comment}

\section{Centraliser Dimension}
\label{section:cd}
\subsection{Definitions and Preliminaries}\label{s:defn}
The centraliser lattices of groups have been studied 
in numerous papers; some 
listed in the introduction. Here we shall consider groups which
have centraliser lattice of finite height; on which there 
is also a considerable literature. We classify such groups according
to centraliser dimension which we define in this section.

If $S$ is a subset of a group $G$ then the centraliser
of $S$ in $G$ is
$C_G(S)=\{g\in G: gs=sg, \textrm{ for all } s\in S\}$. We write
$C(S)$ instead of $C_G(S)$ when the meaning is clear. 
The following properties of centralisers are well-known; see for
example \cite{LR} or \cite{Schmidt94}. Given a family of subsets
$\{S_i\}_{i\in I}$ of $G$ indexed by a set $I$, 
\begin{enumerate}[(i)]
\item\label{clp1} $\cap_{i\in I}C(S_i)=C(\cup_{i\in I} S_i)$;
\item\label{clp2} $\cup_{i\in I}C(S_i)\subseteq C(\cap_{i\in I} S_i)$.
\end{enumerate}
Moreover ({\em loc. cit.}) given subsets $S$ and $T$ of $G$,
\begin{enumerate}[(i)]
\addtocounter{enumi}{2}
\item\label{clp3} if $S\subseteq T$ then $C(S)\supseteq C(T)$;
\item\label{clp4} $S\subseteq C(C(S))$;
\item\label{clp5} $C(S)=C(C(C(S)))$;
\item\label{clp6} $C(S)\subseteq C(T)$ if and only if $C(C(S))\supseteq C(C(T))$.
\end{enumerate}
Let $\fC(G)$ denote the set of centralisers of a group $G$.
The relation of inclusion then defines  a partial order `$\le$' on $\fC(G)$.
We define the infimum of a pair of elements of $\fC(G)$ as
obvious way:
    \[
        C(M_1) \wedge C(M_2)= C(M_1)\cap C(M_2)=C( M_1 \cup M_2).
    \]
Moreover  the supremum $C(M_1) \vee C(M_2)$ of elements
$C(M_1)$ and $C(M_2)$ of $\fC(G)$ may
be
defined to be the intersection of all centralisers containing $C(M_1)$ and $C(M_2)$.
Then  $C(M_1) \vee C(M_2)$ is  minimal among centralisers containing $C(M_1)$ and $C(M_2)$.
These definitions make $\fC(G)$ into a lattice, called the
{\em centraliser lattice} of $G$. This lattice is bounded as it has
a greatest element, $G= C(1)$, and a least element, $Z(G)$, the centre
of $G$. From \eqref{clp1} above, every subset of $\fC(G)$ has an infimum, so
$\fC(G)$ is a complete lattice.

If $C$ and
$C^\prime$
are in $\fC(G)$ with $C$ strictly contained in $C^\prime$ we
write $C<C^\prime$.
If $C_i$ is a centraliser, for $i=0,\ldots ,k$, with $C_0>\cdots >C_k$
then we call $C_0,\ldots, C_k$ a {\em centraliser chain} of {\em
  length} $m$. Infinite descending, ascending and doubly-infinite centraliser
chains are defined in the obvious way. A group $G$ is said to have
the {\em minimal condition on centralisers} min-c if every descending
chain of centralisers is eventually stationary; that is if $\fC(G)$
satisfies the descending chain condition. The {\em maximal
  condition on centralisers} max-c is satisfied by $G$ if 
the ascending chain condition holds in $\fC(G)$. 
From \eqref{clp6} above
a group has min-c if and only if it has max-c.
 Groups with the minimal
condition on centralisers have been widely studied; see
for instance \cite{Zaleski65}, \cite{Schmidt94}, \cite{BryantHartley79}, 
\cite{Bryant79}, \cite{kegel73}. 
As in many of the articles cited we consider now the restriction to
groups in which there is a global bound on the length of centraliser chains.
\begin{defn} \label{def:cd}
If there exists an integer $d$ such that 
the group $G$  has a centraliser chain of length $d$ and no centraliser
 chain of length greater than $d$ then $G$ is said to have {\em centraliser dimension} $\cd (G)=d$.
If no such integer $d$ exists 
we define $\cd (G)= \infty$.
\end{defn}
If $\cd(G)=d$ then every strictly descending chain of centralisers in $\fC(G)$
from $G$
to $Z(G)$ contains at most $d$ inclusions. This number
is usually referred to as the {\em height} of the lattice;
%denoted
%$h(\fC(G))=h(G)$. Thus
so
$\cd (G)$ is the height of the centraliser lattice of $G$.
%=h(\fC(G))$.

Using Definition \ref{def:cd} we introduce the following classes of groups.
For every positive integer $m\ge 0$ set
\[
\DC_m = \left\{ G| \cd(G) \le m \right\}.
\]
In addition we shall sometimes wish to consider the set of all groups with finite
centraliser dimension so we set 
\[
\DC = \bigcup \limits _{i=1} ^\infty \DC_i =\left\{ G| \cd(G) <\infty \right\}.
\]

 Any group from $\DC$ satisfies the minimal condition on
centralisers. The converse is not true: Lennox and Roseblade
\cite[Theorem H]{LR} and Bryant \cite{B} give examples of groups which are nilpotent of
class $2$ and have min-c but are not in $\DC$.
%%%%
The class of groups $\DC$ is nonetheless very broad as the following example shows.
\begin{expl}
\begin{enumerate}[\rm 1.]
  \item Finitely generated Abelian-by-nilpotent groups are in $\DC$ as
    are polycyclic-by-finite groups \cite{LR}.
  \item A linear group of degree $n$ has centraliser dimension at
    most $n^2-1$
    \cite{Weh}. Moreover, if $R$ is a finite direct product of fields then 
    the general linear group $GL(m,R)$ is in $\DC$ \cite{MS}.
  \item If $G$ is a non-Abelian, hyperbolic,
    torsion-free group then, as shown in \cite{MR2}, $G$ is a
    $\CSA$-group so, from Proposition \ref{prop:24}.\ref{prop:clbasic5},
    $\cd(G)=2$.
  \item
    We are grateful to A. Yu. Ol'shanskii for the following argument showing that
    all hyperbolic groups (including those with torsion) have  finite centraliser dimension.
    Suppose $G$ is a hyperbolic group. Then there is a bound on the
    orders of finite subgroups of $G$ (see for example
    \cite[Chapter 4]{GH}). Thus it suffices to show that there is a
    bound on the length of strictly descending chains of infinite
    centralisers of $G$. Suppose that $C=C_G(X)$ where $X$ is a subset
    of $G$ generating a non--elementary subgroup $K$ of $G$. 
    (A group is {\em elementary} if it has a cyclic subgroup of finite
    index.
    Also, the elementariser $E_G(H)$ of a subgroup $H$ of $G$ is the set of all
    $x\in G$ such that $x^H$ is finite: see \cite{IO1} for details.)
    Then, from  \cite[Proposition 1]{IO1},
    the elementariser $E_G(K)$ of $K$ is finite, so $C=C_G(K)\subseteq
    E_G(K)$ is finite.
    Therefore infinite centralisers in $G$ are centralisers
    of elementary subgroups of $G$. Since there is a bound on the
    order of finite subgroups of $G$ the set of lengths of chains of
    centralisers of finite elementary subgroups is bounded. As
    $E\subseteq E^\prime$ implies $C(E)\supseteq C(E^\prime)$ it
    therefore suffices to bound the set of lengths of chains of infinite
    centralisers of infinite elementary subgroups. Suppose now that 
    $C_G(E)$ is the centraliser of the elementary subgroup $E$ and
    that $a$ is an element of infinite order in $E$. Then $C_G(a)$ is 
    an (infinite) elementary subgroup (\cite[p. 156]{GH}) and $C_G(a)\supseteq C_G(E)$, so
    $C_G(E)$ is elementary. Let 
    \begin{equation}\label{eq:burnc}
    C_0\sdc \cdots \sdc C_d
    \end{equation}
    be a strictly descending chain of infinite centralisers
    $C_i=C_G(E_i)$, where $E_i$ is an infinite elementary subgroup. We
    may assume that $E_0<\cdots <E_d$. There is an element of infinite
    order in $C_0\cap E_d$ so the group $E$ generated by $C_0$ and $E_d$
    is elementary (\cite[p. 375]{IO1}). Thus \eqref{eq:burnc} is a
    chain of centralisers of the elementary subgroup $E$. From
    \cite[Lemma 19]{IO1} it follows that there is an integer $M$ such
    that every infinite elementary subgroup of $G$ has an infinite cyclic 
    subgroup of index at most $M$. From this and Proposition
    \ref{prop:ftext} below it follows  that $G$ is in $\DC$, as claimed.
  \item We are grateful to S. V. Ivanov for the following argument
    concerning free Burnside groups.
    The free Burnside groups of large exponent have
    centraliser dimension $2$, when $n$ is odd, but do not have
    min-c when $n$ is even. 
    In more detail, let
    $G=B(m,n)$ be the $m$-generator free Burnside group of exponent $n$. 
    If $m>1$ and  $n\geq 665$ then centralisers of non-trivial
    elements of $G$ are cyclic of order $n$ \cite{A}. It follows that
    in this case $\cd{G}=2$. On the other hand suppose that $n\ge
    2^{48}$ and that $2^9|n$. Then we may choose a finite $2$-subgroup
    $T_1$ of $G$. From \cite[Theorem 1 (a)]{IO2}, $C_G(T_1)$ contains
    a subgroup $B$ isomorphic to $B(2,n)$ such that $C_G(T_1)\cap
    B=\{1\}$. We may now take a finite $2$-subgroup $D$ of $B$ and
    set $T_2=\langle T_1,D\rangle =T_1\times D$. Then $T_2$ is a
    finite $2$-subgroup of $G$. Repeating the process starting with
    $T_2$ instead of $T_1$ and continuing this way we see that $G$
    contains an infinite ascending chain 
    \[T_1<T_2<\cdots \]
    of finite $2$-subgroups. From \cite[Theorem 1 (c)]{IO2}, for all
    $i$, $C_G(C_G(T_i)=T_i$, so 
    \[C_G(T_1)>C_G(T_2)>\cdots \]
    is an infinite descending chain of centralisers. Therefore $G$
    does not have min-c.
    \item We pose the following question.
      Is the centraliser dimension of a biautomatic group finite?
      This is related to (and stronger than) several well-known questions
      concerning these groups.
      Gersten and Short \cite[Proposition 4.3]{GS} show that,
      in a biautomatic group, centralisers of finite subsets are
      biautomatic. They ask (loc. cit.) if biautomatic
      groups have min-c and show that if so then every Abelian subgroup of a
      biautomatic group is finitely generated. Moreover Mosher
      \cite{M} shows that a biautomatic group has an infinitely generated
      Abelian subgroup if and only if it has an Abelian subgroup which
      is either of infinite rank or is an infinite torsion
      group. However, whether or not such subgroups are to be found in
      biautomatic groups is an open question. Another related
      open question asks whether or not a biautomatic group can have an element
      of infinite order which has infinite index in its
      centraliser. Mosher (loc. cit.) shows that if such a biautomatic group
      exists it must contain a subgroup isomorphic to $\ZZ^2$, so the
      group cannot be hyperbolic.
    \item Blatherwick \cite{B} has examples showing that, for each integer
      $m\ge 4$ (and for $m=2$) there exists a nilpotent
      group of class $2$ with
      centraliser dimension $m$  (see also \cite{LR}).
    \item We shall show in Section \ref{subs:behave} that the class
      $\DC$ is closed under formation of  direct sums and free products
      (with finitely many factors) and certain amalgamated
      products. Moreover if a group $G$ has a subgroup of finite index
      belonging to $\DC$ then $G$ is in $\DC$ (Proposition \ref{prop:ftext}).
 \end{enumerate}
\end{expl}
\begin{comment}
\begin{con}
We conjecture that biautomatic groups
have finite centraliser
    dimension.
\end{con}
\end{comment}
%%%%
The first four statements of the following proposition are well-known,
but we give proofs for completeness. Statements \ref{prop:clbasic1}
and \ref{prop:clbasic2} follow from Lemma 3.1 and Theorem 3.2 of
\cite{Schmidt94}, which show that a group has distributive centraliser
lattice if and only if the group is Abelian, in which case the lattice
is trivial; and that no group can have centraliser lattice of height one.
\begin{prop}\label{prop:cdbasic}
\begin{enumerate}[\rm 1.]
\item\label{prop:cdbasic1}  
If $G$ has min-c and 
$C$ is a centraliser in $G$
then there exists a finite
subset $M$ such that $C=C(M)$ \cite{Bryant79}.
\item\label{prop:cdbasic2}  If $\cd(G)=m$ and
\[
G=C_0\sdc \cdots \sdc C_{m}=Z(G)
\]
is a centraliser chain of maximal length in $G$ then
$C_{m-1}$ is Abelian \cite{Schmidt94}.
\item\label{prop:clbasic1} Let $G$ be an Abelian group, then
  $\cd(G)=0$ \cite{Schmidt94}.
    \item\label{prop:clbasic2} If $G$ is non-Abelian then $\cd(G) \ge
      2$: that is $\DC_0=\DC_1$ \cite{Schmidt94}.
    \item\label{prop:clbasic4} In the event that $\cd(G)=m$ is finite, there exists
    an $m$-tuple of non-central elements $a_1,\ldots ,a_{m}$ such that
    \[
 G  \sdc C(a_1) \sdc \ldots \sdc C(a_1, \ldots, a_{m})=Z(G).
    \]
\end{enumerate}
\end{prop}
\begin{proof} If $C=C(S)$ is not the centraliser of any finite subset then we
may construct an infinite centraliser chain by choosing succesive elements $s_1,s_2,\ldots $ of
$S$ and forming centralisers $C(s_1,\ldots, s_k)$, for increasing
$k$.
This proves \ref{prop:cdbasic1}. To see \ref{prop:cdbasic2}
suppose that $C_{m-1}$ is
non-Abelian. Take a
pair $a,b$ of non-commuting elements from $C_{m-1}$
and consider the centraliser $C=C_{m-1}\cap C(a)$.
Notice that $a\in C$ but, since $[a,b]\ne 1$, $b\notin C$.
Hence $C_{m-1}\sdc C$ and we have a centraliser chain
\[
G=C_0\sdc \cdots \sdc C_{m-1}\sdc C\sdc C_{m}=Z(G)
\]
of length greater than $m$. As $\cd(G)=m$ this is a contradiction and
\ref{prop:cdbasic2} holds.

%%%
Statement \ref{prop:clbasic1} is clear.
For \ref{prop:clbasic2} observe that if $G$ is non-Abelian then
$G\neq Z(G)$ and so we may choose $a\in G\backslash Z(G)$. Then
\[
G  \sdc C(a) \sdc Z(G),
\]
so $\cd (G) \ge 2$.

To prove statement \ref{prop:clbasic4} notice that there is nothing to prove if $m<2$. Assume
that $\cd(G) = m\ge 2$. Then there exist
finite subsets $M_1,\ldots , M_{m}$ of $G$ such that
\begin{equation}\label{e:sdc}
G\sdc C(M_1)  \sdc C(M_2)  \sdc \ldots \sdc C(M_{m})=Z(G)
\end{equation}
is a (strictly descending) centraliser chain.
Take  $a_1 \in M_1$ such that $C(a_1) \ne G$. Then
$C(a_1)=C(M_1)$ by maximality of \eqref{e:sdc}. Assume that
elements
$a_1,\ldots ,a_{i-1}$ have been chosen so that $C(M_j)=C(a_1,\ldots, a_j)$,
for $j=1,\ldots ,i-1$. As $C(M_i)\sac C(M_{i-1})$ we may choose $a_i \in
M_i$ such that $C(a_i)\ngeq C(M_{i-1})$.  Since $C(a_i)\sdc C(M_i)$ and
\eqref{e:sdc} is maximal it follows that $C(M_i)=C(a_1,\ldots, a_i)$.
Hence, by induction,  we may choose such $a_i$ for $i=1,\ldots ,m$. None of
the $a_i$ belong to $Z(G)$, since \eqref{e:sdc} is a strictly
descending chain, hence \ref{prop:clbasic4} 
holds.
\end{proof}

From now on we shall only consider groups with finite centraliser
dimension.  
Next we show that for every
positive integer $m\ge 0$ the class of groups $\DC_m$  is
universally axiomatisable. We shall make use of the notation
of Section \ref{s:logic} for formulas in the language $L$.
Since $\DC_0=\DC_1$ and $\DC_0$ is the class of all Abelian
groups, these classes are defined by the following universal
sentence,
\[
\DC_0=\DC_1: \forall  x,y ( \left[ x,y \right]=1),
\]
which, in the notation of Section \ref{s:logic}, takes the form
\[
\forall x,y (x \in C(y)).
\]
We next write down an axiom for $m=2$.
\begin{multline*}
  \DC_2: \forall x_0,x_1,x_2,y_1,y_2,z\
  (y_1 \in C(x_0) \wedge  y_1 \notin C(x_0,x_1) \wedge \\
  \wedge y_2 \in C(x_0,x_1) \wedge y_2
  \notin C(x_0,x_1,x_2)\wedge z \in
  C(x_0,x_1,x_2) %\rightarrow \\
  \rightarrow  z\in Z(G)).
\end{multline*}
From Proposition \ref{prop:cdbasic}.\ref{prop:clbasic4} it follows that $\DC_2$ is the class of groups
axiomatised by this universal sentence.
For $m>2$ we have the following axiom.
\begin{gather*}
\DC_m  : \forall x_0,\ldots,x_m,y_1,\ldots,y_{m},z\ ( y_1 \in C(x_0) \wedge y_1 \notin C(x_0,x_1)\wedge\\
\wedge y_2 \in C(x_0,x_1) \wedge y_2  \notin C(x_0,x_1,x_2) \wedge\ldots \wedge \\
\wedge y_{m} \in C(x_0,\ldots, x_{m-1})\wedge y_{m} \notin C(x_0,\ldots, x_m)\\
%\wedge%\\
\wedge z \in C(x_0,\ldots, x_m) \rightarrow z\in Z(G)).
\end{gather*}
\begin{prop}
For $m\ge 0$, the class of groups  $\DC_m$,
defined in Section \ref{s:defn}, is axiomatised by the universal sentence
$\DC_m$ above. That is, a group $G$ satisfies axiom $\DC_m$
if and only if
 $\cd (G)\le m$.
\end{prop}
\begin{proof}
Obviously, if $\cd(G) \le m$ then $\DC_m$ holds in $G$.
Conversely, suppose that $\DC_m$ holds in $G$ and $\cd(G) = n>m$. Then,
by Proposition \ref{prop:cdbasic}.\ref{prop:clbasic4} there exist non-central
elements $a_1, \ldots, a_{n}$ such that
\[
G\sdc
C(a_1)\sdc C(a_1, a_2) \sdc \ldots \sdc C(a_1, \ldots,a_{n})=Z(G).
\]
Therefore there are elements $y_1,\ldots ,y_{n-1}\in G$ such that
$y_1\in G$, $y_1\notin C(a_1)$ and $y_j\in C(a_1,\ldots ,a_{j-1})$, $y_j\notin C(a_1,\ldots ,a_j)$,
for $j=2,\ldots ,n$. In this case though
$\DC_m$ does not hold for the $(2m+1)$-tuple $1,a_1,\ldots,
a_{m},y_1,\ldots y_m$ of elements of $G$.
\end{proof}

\subsection{The Behavior of Centraliser Dimension Under Group
Operations}\label{subs:behave}

In the proof of Proposition \ref{prop:ofcd} we shall make use of
the following well-known theorem.

\begin{thm}[Theorem 4.5 of \cite{MKS}] \label{thm:comcr}
Let $G$ be an amalgamated product of
$H_1$ and $H_2$ with amalgamation by $K$, i.e.
$G=  H_1 \ast_{K} H_2$.
Assume that the elements $q$ and $r$ of $G$ commute. Then one of the following holds.
\begin{enumerate}[\rm 1.]
\item Either $q$ or $r$ is conjugate to an element from the
  subgroup $K$.
\item Both $q$ and $r$ lie
  in $H_i^{g}$, for some $g\in G$ and some
  fixed $i=1$ or $2$.
\item If the previous two conditions fail then
  \[
  q= g c g^{-1} u^l \; \hbox{ and } \; r= g c' g^{-1} u^{l'}
  \]
  where $c,c' \in K$, $g, u \in G$, $l,l' \in
  \ZZ$ and the elements  $gcg^{-1}$, $gc'g^{-1}$ and $u$
  pairwise commute.
\end{enumerate}
\end{thm}

In order to state our
next proposition we need a further definition. Let $Z(G)$ denote
the centre of the group $G$. We define
\[
\ttz(G)= \left\{
  \begin{array}{ll}
    0, & \textrm{if } Z(G)=1_G\\
    2, & \textrm{otherwise}
  \end{array}
\right. .
\]
Statement \ref{i:cd2} of the following proposition 
is a consequence of the fact, proved in \cite{Schmidt95}, that 
$\fC(G_1\times G_2)=\fC(G_1)\times \fC(G_2)$,
for
groups $G_1$ and $G_2$. We give a proof for completeness.
\begin{prop}\label{prop:ofcd} Let $G_1,G_2 \in \DC$.
~ \be[\rm 1.]
    \item\label{i:cd1} If $G_1 \le G_2$
    then $\cd(G_1) \le \cd(G_2)$.

    \item\label{i:cd2} $\cd(G_1 \times G_2)=\cd(G_1)
      + \cd(G_2)$ (see \cite{Schmidt95}). 
    \item\label{i:cd3}
$
    \cd(G_1 \ast G_2)=
    \max\left\{ \cd(G_1)+\ttz(G_1),\cd(G_2)+\ttz(G_2)\right\}.
$ \item\label{i:cd4} Let $G_1$ and $G_2$  be non-Abelian groups
with $Z(G_i)=Z_i$, $i=1,2$, such that $Z_1\simeq Z_2\ne 1$, and
let $G=G_1 \ast_{Z_1=Z_2} G_2$. Then
\[
    \cd(G)=
    \max\left\{ \cd(G_1),\cd(G_2)\right\}.
\]
\ee
\end{prop}
\begin{proof} \ref{i:cd1}. This follows from the fact that every
 centraliser $C_{G_1}(M)$ of a
set $M$ in $G_1$ is a subset of the centraliser $C_{G_2}(M)$ in
$G_2$. Moreover if $C_{G_1}(M) \sac  C_{G_1}(N)$ then
$C_{G_2}(M) \sac  C_{G_2}(N)$.\\[1em]
\ref{i:cd2}.
Let $G=G_1\times G_2$. If $M\subseteq G$ and we let $M_1$ and
$M_2$ be the projections of $M$ onto $G_1$ and $G_2$, respectively,
then $C_G(M)=C_{G_1}(M_1)\times C_{G_2}(M_2)$. In particular,
if $N_i\subseteq G_i$, for $i=1,2$, then $C_G(N_1\times N_2)=C_{G_1}(N_1)\times C_{G_2}(N_2)$.

We first show that $\cd(G) \ge \cd(G_1)
+ \cd(G_2)$. Let $\cd(G_1)=m$ and $\cd(G_2)=n$. In this case there exist
centraliser
chains
\[
G_1\sdc C_1\sdc \cdots \sdc C_{m}=Z(G_1)
\]
and
\[
G_2\sdc D_1\sdc \cdots \sdc D_{n}=Z(G_2).
\]
Then
\[
G\sdc C_1\times G_2\sdc \cdots \sdc C_{m}\times G_2 \sdc
C_{m}\times D_1\sdc \cdots \sdc C_{m}\times D_{n}
\]
is a strictly descending chain of centralisers in
$G$.
Therefore $\cd(G) \ge \cd(G_1) + \cd(G_2) $.

Next we show that $\cd(G) = \cd(G_1) + \cd(G_2) $.
Suppose we have a centraliser chain
\begin{equation}\label{e:dpchain}
C_0\sdc \ldots \sdc  C_k
\end{equation}
 in $G$.
Then $C_i = C_{G_1}(M_i) \times C_{G_2}(N_i)$,
where $M_i \subseteq G_1$ and $N_i \subseteq G_2$, for $i=0 ,\ldots ,k$.
Since \eqref{e:dpchain}
is strictly descending we have
$C_{G_1}(M_i) \edc C_{G_1}(M_{i+1})$ and  $C_{G_2}(N_i) \edc C_{G_2}(N_{i+1})$,
with at least one of these inclusions strict, for $i=0,\ldots ,k-1$.
Since $\cd(G_1)=m$ there are at most $m+1$ distinct centralisers among the
$C_{G_1}(M_i)$. Hence there are at most $m$  of these inclusions with
$C_{G_1}(M_i) \ne C_{G_1}(M_{i+1})$. Similarly there are at most $n$ inclusions
with $C_{G_2}(N_i) \ne C_{G_2}(N_{i+1})$. Hence  the number $k$ of inclusions in \eqref{e:dpchain}
is at most $m+n$, and it follows that $\cd(G)\le m+n$.\\[1em]
\ref{i:cd3}. Let $G=G_1*G_2$
\begin{comment}
 and suppose first that
$G_1$ and $G_2$ are Abelian. Then $G$ is a
$\CT$-group with trivial centre and, by Proposition \ref{prop:24}.\ref{prop:clbasic5}, $\cd(G)=2$.
Therefore the statement holds if $G_1$ and $G_2$ are Abelian.
We may therefore assume that at least one of $G_1$ and $G_2$
is non-Abelian. Let 
\end{comment}
and let 
$f$ and  $g$ be two non-trivial elements of
$G$ such that $C(f)\ne C(g)$. Then either $C(f) \cap C(g) = 1$ or
both of these elements and their centralisers lie in
$G_i^h$, for some fixed $h\in G$ (see Theorem \ref{thm:comcr}). This implies that if
\begin{equation}\label{eq:ac}
G \sdc C_1 \sdc  \cdots \sdc C_p \sdc 1
\end{equation}
is a
strictly descending chain of centralisers with $p\ge 2$
then there are  fixed $i$ and $h$ such that $C_j \eac G_i^h$, for all
$j$. After conjugation by $h^{-1}$ we may then assume that $C_j\eac G_i$, for all $j$.
(If $p=1$ then we may replace $C_1$ with the centraliser of an
element of $G_1$ or $G_2$, if necessary; so we may assume
that the claim holds in this case as well.)

First suppose that 
$G_i$ has trivial centre. Then, for $a\in G_i$,
$C_{G}(a) \ge  G_i$ only if  $a=1$ in which case $C_G(a)=G$.
Hence 
$C_1$ is strictly contained in  $G_i$. Also, $C_p\neq 1$ implies $C_p\neq Z(G_i)$.  
Therefore
\begin{equation*}%\label{eq:Gichain}
G_i \sdc C_1 \sdc  \cdots \sdc C_p \sdc 1
\end{equation*}
is a centraliser chain of length $p+1$ in $G_i$ and we have 
$p+1\le \cd(G_i)= \cd(G_i)+\ttz(G_i)$.

On the other hand, if $Z(G_i)\ne 1$
then $G_i=C_G(g)$, for all $1 \ne g \in Z(G_i)$. In this case it is
possible that $C_1=G_i^h$ and $C_{p}=Z^h(G_i)$. Thus, if $Z(G_i)\ne 1$ then
$p+1\le \cd(G_i)+2= \cd(G_i)+\ttz(G_i)$.
It follows that $\cd(G)\le \max\left\{ \cd(G_1)+\ttz(G_1),\cd(G_2)+\ttz(G_2)\right\}.$

Conversely suppose that 
\[
G_i\sdc C_1\sdc \cdots \sdc C_p=Z(G_i)
\]
is a centraliser chain in $G_i$. If $Z(G_i)=1$ then, replacing $G_i$ with
$G$ in this chain we obtain a centraliser chain for $G$ of length $p$. 
If $Z(G_i)\neq 1$ then adding $G$ to the left and $1$ to the right of this chain
we obtain a centraliser chain for $G$ of length $p+2$. Hence
$\max\left\{ \cd(G_1)+\ttz(G_1),\cd(G_2)+\ttz(G_2)\right\}\le \cd(G)$.\\[1em]
\ref{i:cd4}.  We have $Z(G)=Z(G_1)=Z(G_2)$ and  Theorem \ref{thm:comcr} takes the following
form. If $x,y$ are elements from $G$ such that $xy=yx$
then
\begin{enumerate}[\rm (i)]
    \item\label{i:comm1} $x$ or $y \in Z(G)$; or
    \item\label{i:comm2} there is $g\in G$ such that $x \in
      G_i^g\smallsetminus Z(G)$ and  $y\in G_i^g\backslash Z(G)$; or
    \item (\ref{i:comm1}) and (\ref{i:comm2})
      do not hold, and there exists
    an element $z$, such that $x=c_1z^k$ and $y=c_1z^l$, with $c_1,c_2
    \in Z(G)$, $k,l \in \ZZ$.
\end{enumerate}

Now let $f, g$ be two elements from $G\smallsetminus Z(G)$ such
that $C_G(f)\ne C_G(g)$. Then either $C_G(f)\cap C_G(g)= Z(G)$ or $f$
and $g$ both lie in the same subgroup of the type $G_i^h$. Since
neither of $G_1$ or $G_2$ is a centraliser in $G$
the result follows as in \ref{i:cd3}.
  \end{proof}

As shown in \cite{Schmidt94},
if $G=HK$ where $H\cap K\neq \nul$ then it is not necessarily the 
case that $\fC(G)\cong \fC(G)\times \fC(G)$ even if $H$ and $G$
centralise
one another.
Also the relationship between $\fC(G)$ and $\fC(G/Z(G))$ is
complicated: \cite{Bryant79} contains an example of a group $G$ such
that
$G$ has min-c  but $G/Z(G)$ does not. Moreover Example
\ref{ex:nilp}.\ref{ex:nilp3}
below shows that centraliser dimension may increase on factoring 
by the centre.
On the positive side it is shown in \cite{Schmidt94} that if
$G=HZ$, where $Z\le Z(G)$, then $\fC(G)\cong\fC(H)$; so
$\cd(G)=\cd(H)$. We also have the following proposition.
\begin{prop}\label{prop:nilab}
Let $G$ be a group such that every nilpotent subgroup of $G$ is
Abelian. Then $\cd(G)=\cd(G/Z(G))$.
\end{prop}
\begin{proof}
 Let  $\bar G$ denote $G/ Z(G)$ and
$\bar g$ the image of $g$ in $\bar G$. If  $g\in G$ then $C(\bar
g)=\left\{\bar f| \left[ g,f \right]\in Z(G) \right\}$. Let $g,
f\in G$ such that $\bar f\in C(\bar g)$. Let $H=\langle
f,g\rangle$. Then $[[f,g],h]=1$, for all $h\in H$, so $H$ is
nilpotent. Hence $H$ is Abelian, so $f\in C(g)$. Therefore $C(\bar
g)=C(g)/Z(G)$ and it follows, via a straightforward induction, that
$C(\bar X)=C(X)/Z(G)$, for all finite sets $X\subseteq G$. Therefore
 $\cd(\bar G)=\cd(G)$.
\end{proof}

The class of groups
satisfying min-c
is closed under the formation 
of finite extensions \cite{kegel73}. The same is true of  the class $\DC$, however 
\begin{prop}\label{prop:ftext}
in this case slightly more can be said.
Let $H$ be a subgroup of finite index $k$ in a group $G$. If
$\cd(H)=d<\infty$ then $\cd(G)\le ((d+2)k+2)k$
\end{prop}
\begin{proof}
Let 
\[G=C_0>C_1>\cdots >C_n=Z(G)\]
be a centraliser chain in $G$ of length $n$
and let elements $t_1,\ldots ,t_k$ of $G$ form a transversal for $H$
in $G$, with $t_1=1$. Then 
\begin{displaymath}
C_j=\bigcup_{i=1}^k(C_j\cap Ht_i), \textrm{ for } j=0,\ldots , n.
\end{displaymath}
%Since $C_j\cap Ht_i=\emptyset$ implies $C_l\cap Ht_i=\emptyset$, for
%all $l\ge j$, we may define 
Note that if $C_j\cap Ht_i=\emptyset$ then  $C_l\cap Ht_i=\emptyset$,
for
all $l\ge j$.

For $i=1,\ldots ,k$ define $d(i)=n$, if $C_n\cap Ht_i\neq \emptyset$,
and  otherwise $d(i)=j$, where $j$ is the unique integer such that $C_j\cap
Ht_i\neq \emptyset$ and $C_{j+1}\cap Ht_i=\emptyset$.
Since $t_1=1$ and $H\cap Z(G)\neq \emptyset$ we have $d(1)=n$. We may
therefore reorder the $t_i$'s so that \[n=d(1)\ge d(2)\ge \cdots \ge
d(k)\ge 0.\] 
Now fix $j$ with $0\le j<n$. If $j=d(s)$, for some $s$, then 
\[C_j=\bigcup_{i=1}^s(C_j\cap Ht_i) \textrm{ and
}C_{j+1}=\bigcup_{i=1}^r(C_{j+1}\cap Ht_i),\]
where $r<s$ is maximal such that $d(r)\ge j+1$.
On the other hand, if $j\neq d(s)$, for all $s$, then there exists $s$
such that $d(s)>j>d(s+1)$. Then, since $C_j\cap Ht_{s+1}=\emptyset$ by
definition of $d(s+1)$,
\[C_j=\bigcup_{i=1}^s(C_j\cap Ht_i) \textrm{ and
}C_{j+1}=\bigcup_{i=1}^s(C_{j+1}\cap Ht_i).\]
Moreover, by definition of $d(s)$, $C_{j+1}\cap Ht_s\neq \emptyset$.
As $C_j>C_{j+1}$ it follows that, for some $i$ with $1\le i\le s$, 
\begin{equation}\label{eq:cosC}
C_j\cap Ht_i>C_{j+1}\cap Ht_i.
\end{equation}
For $j=0,\ldots,n$ define $e(j)=0$, if $j=d(i)$, for some $i$, and
otherwise set
$e(j)=i$, where $i$ is chosen to satisfy \eqref{eq:cosC}. Then
$|e^{-1}(0)|\le k$ and $\sum_{i=0}^k|e^{-1}(i)|=n$.  If   $l$ is an
integer such that  $n\ge (l+1)k+1$ this implies that there is some
$s$, with $1\le s\le k$, such that $|e^{-1}(s)|\ge l$.
Assume then that $n\ge (l+1)k+1$, for some positive integer $l$ and
fix such an $s$. Then there are integers $j_1\le \cdots\le j_l$ such
that  
$e(j_r)=s$, for $r=1,\ldots ,l$. From \eqref{eq:cosC} it follows that 
\begin{equation}\label{eq:cosB}
C_{j_1}\cap Ht_s\sdc \cdots \sdc C_{j_l}\cap Ht_s.
\end{equation}
Now $ C_{j_l}\cap Ht_s\neq \emptyset$, so there exists an element $y\in
 C_{j_l}\cap Ht_s$. As $y\in C_{j_r}$ and $Ht_s=Hy$ we have 
\[
C_{j_r}\cap Ht_s=C_{j_r}\cap Hy=(C_{j_r}\cap H)y,
\]
for $r=1,\ldots ,l$. Hence \eqref{eq:cosB} implies that 
\begin{equation}\label{eq:cosA}
C_{j_1}\cap H\sdc \cdots \sdc C_{j_l}\cap H.
\end{equation}

%%%%%%%%%%%%%%
For
$r=1,\ldots ,l$ we have $C_{j_r}=C_G(X_r)$, where $X_r\subseteq G$. 
We may assume that $X_1\subseteq X_2\subseteq \cdots \subseteq X_{l}$ 
and set $Y_r=X_{r-1}\backslash X_r$, for $r=2,\ldots l$.
Next we shall argue that we may assume that $Y_r$ has only one element, for all $r$.
To see this suppose that $r$ is minimal such that $|Y_r|>1$. There exists 
$h\in C_G(X_{r-1})\cap H$ such that $h\notin C_G(X_{r})\cap H$. As $X_r=X_{r-1}\cup Y_r$
this means that there is $y\in Y_r$ such that $h\notin C_G(y)$. Hence
$C_G(X_{r-1})\cap H> C_G(X_{r-1})\cap C_G(y)\cap H\ge C_G(X_{r})\cap H$. 
Thus we may replace $Y_r$ with $\{y\}$. Continuing this way each $Y_r$ may be replaced
by a singleton. 
We may now assume that there are elements $x_1,\ldots ,x_l$ of $G$ 
such that 
$X_r=\{x_1,\ldots ,x_r\}$, for
$r=1,\ldots ,l$. 

Define $c(r)=i$ to be the unique integer such that
$x_r\in Ht_i$.
Let $m$ be an integer such that $l\ge mk+1$. Then for some $i$, with
$1\le i\le k$, we
have $|c^{-1}(i)|\ge m$. Fix $s$ such that $|c^{-1}(s)|\ge m$ and let $c^{-1}(s)\supseteq
\{r_1,\ldots ,r_m\}$, where $r_1<\cdots <r_m$. Set
$x_{r_i}=y_i$ and $Y_i=\{y_1,\ldots y_i\}$, for $i=1,\ldots ,m$. Then it follows that
\begin{equation}\label{eq:ydown}
C_G(Y_1)>\cdots >C_G(Y_m)
\end{equation}
is a strictly descending chain of centralisers in $G$. Since $c(r_i)=s$
we have $y_i=a_it_s$, where $a_i\in H$, for all $i$. A straightforward calculation
shows that, for any elements $a, b,c\in G$, 
the identity $C_G(ab)\cap C_G(cb)=C_G(ab)\cap C_G(ac^{-1})$ holds. By induction it
follows that \[C_G(Y_i)=C_G(a_1t_s)\cap C_G(a_1a_2^{-1},\ldots ,a_1a_i^{-1}),\]
for $i=2,\ldots ,m$. Setting $A_i=\{a_1a_2^{-1},\ldots ,a_1a_i^{-1}\}\subseteq H$, 
it follows 
from \eqref{eq:ydown} that
\[G>C_G(A_2)>\cdots >C_G(A_m)\]
is a strictly descending chain of centralisers in $G$. 

Now define $A_1=\{1\}$ and $D_i=C_G(A_i)$, for $i=1,\ldots ,m$.
Then $D_i\cap H=C_H(A_i)$, so 
\[H=H\cap D_1\ge \cdots \ge H\cap D_m\]
is a strictly descending chain of centralisers of length $m-1$ in $H$.
This occurs if $n\ge (l+1)k+1$ and $l\ge mk+1$; that is $n\ge
(mk+2)k+1$. Thus if $n\ge ((d+2)k+2)k+1$ we obtain a contradiction,
and the result follows.
\end{proof}

\subsection{Groups of Centraliser Dimension $2$}\label{subs:dim2}

In this section we concentrate attention on
the class $\DC_2$ of groups that have centraliser dimension at most
$2$. There are many examples of such groups:
free groups, torsion-free hyperbolic 
groups and free Burnside groups, 
of large odd exponent, have centraliser dimension $2$. R.~Schmidt
\cite{Schmidt94} has completely classified finite groups of 
centraliser dimension $2$ (which are ${\mathfrak{M}}$-groups
in the terminology of \cite{Schmidt94}).
Locally finite groups in
$\cd(G)=2$ have also been fairly intensively studied
(see Chapter 2 of \cite{book}, and 
\cite{busgor21}). Here  we show that there is a connection between
groups with centraliser dimension $2$ and $\CT$-groups and give some 
examples.

\begin{prop} \label{prop:24}
Let $G$ be a non-Abelian group.
\begin{enumerate}[\rm 1.]
\item\label{prop:clbasic5} If $\cd(G)=2$ and $Z(G)=1$ then $G$ is a 
  $\CT$-group. Conversely,
  if $G$ is a $\CT$-group then  $\cd(G)=2$ and $Z(G)=1$.
\item\label{prop:2}  Suppose that every nilpotent subgroup of $G$ is
  Abelian. Then $\cd(G)=2$ if and and only if the factor-group $G/Z(G)$ is a $\CT$-group.
\end{enumerate}
\end{prop}
\begin{proof}
To see \ref{prop:clbasic5} suppose there exists a non-Abelian
$\CT$-group $G$ such that $\cd(G) \ge 3$. Then, from Proposition
\ref{prop:cdbasic}.\ref{prop:clbasic4}, there exists a chain of centralisers
\[
G  \sdc C(a_1) \sdc C(a_1, a_2) \sdc Z(G).
\]
Since the second inclusion above is strict it follows that
$C(a_1)\ne C(a_2)$, and since $G$ is a $\CT$-group this implies
$[a_1,a_2]\ne 1$. As the third inclusion is strict there is an non-trivial
element  
$b \in C(a_1) \cap C(a_2)$. The assumption
that $b \ne 1$ together with the $\CT$ axiom now imply that
$\left[a_1, a_2 \right]= 1$, a contradiction.

To prove the converse suppose that $\cd(G)=2$ and that $G$ is
non-Abelian group with trivial centre. Since $G$ is non-Abelian
the centraliser of a non-trivial element is a proper non-trivial
subgroup. As $\cd(G)=2$ such centralisers are all Abelian
subgroups, by Proposition \ref{prop:cdbasic}.\ref{prop:cdbasic2}.
Now if $b_1$ and $b_2$ belong to the centraliser of a non-trivial
element $a \in G$ then $b_1$ and $b_2$ commute, so the $\CT$ axiom
is satisfied.

%Since $G$ has no non-Abelian nilpotent subgroups the centre
%of $\bar G= G/ Z(G)$ is trivial.
%%
In the setting of \ref{prop:2} note that, from 
Proposition \ref{prop:nilab}, $\cd(G)=\cd(G/Z(G))$ and also that
$G/ Z(G)$ is non-Abelian, for
otherwise $G$ is nilpotent and thus, 
by hypothesis Abelian. Suppose first that $\cd(G)=2$.  
Since $G$ has no non-Abelian nilpotent subgroups the argument of the 
proof of Proposition \ref{prop:nilab} shows that the centre
of $G/ Z(G)$ is trivial.
Therefore \ref{prop:clbasic5} implies that $G/ Z(G)$ is a
$\CT$-group. On the other hand, if $G/ Z(G)$ is a
$\CT$-group then it follows, from \ref{prop:clbasic5}, 
that $\cd(G)=2$.
\end{proof}
\begin{rem}
In the event that $\cd(G)=2$ and $Z(G) \ne 1$ the situation
is much more complex (see Example \ref{ex:nilp}.\ref{ex:nilp3} below).
\end{rem}

\begin{expl}\label{ex:nilp}
\begin{enumerate}[\rm 1.]
    \item Let $G$ be a non-Abelian $\CT$-group and $A$ be an Abelian group
    then, by Propositions \ref{prop:24} and  \ref{prop:ofcd}, $\cd(G\times A)=2$.
    \item For every positive integer $c\ge 3$ there exists a
    nilpotent group  $G$ of class $c$ such that
    $\cd(G)=2$. Here is an example of such group. Let
    $A=\ZZ^c$ be a lattice of the rank $c, c \ge 2$. Take
    an automorphism of $\phi$ of $A$ given by the following matrix of order
    $c$ in the natural base $e_1, \ldots, e_c$:
\[
\left[ \phi \right]=
    \left(%
\begin{array}{ccccc}
  1 & 1 & \cdots & 0 & 0 \\
  0 & 1 & \cdots & 0 & 0 \\
    &   & \cdots &   &   \\
  0 & 0 & \cdots & 1 & 1 \\
  0 & 0 & \cdots & 0 & 1 \\
\end{array}%
\right)
\]
Let $\a$ be the homomorphism from the infinite cyclic group $\la t
\ra$ to Aut$(A)$ given by $\a(t)=\phi$ and
let $G$ be a semi-direct product %of $A$ by $\phi$: that is
$G=\la t \ra {\ltimes}_\a A$. Then $Z(G)=\left<e_c \right>$ and
the upper central series of $G$ is
\[
1\le \la e_c \ra \le \la e_{c-1},e_c\ra \le \cdots \le \la
e_2,\cdots,e_c\ra \le G,
\]
so $G$ is a nilpotent group of class $c$. We claim that centralisers
in $G$ are either $G$, $A$, $Z(G)$ or of the form $\la t^ia,
e_c\ra$, for some $i\in \ZZ$ and $a\in A$. To see this first note
if $g\in A$ and $g\notin \la e_c\ra$ then the centraliser of $g$
is $A$. Moreover if $g\in \la t\ra$ then the centraliser of $g$ is
$\la t,e_c\ra$. Hence it remains to calculate the centraliser of
$t^r a$, where $0\neq r\in\ZZ$ and $a\in A$, $a\notin\la e_c\ra$.
Note that $G$ is torsion-free so extraction of roots in $G$ is
unique. (For $h\in G$ and $n\in \NN$, the equation $x^n=h$ has at
most one solution. See for example \cite{K}.) Hence, for $x,c\in
G$, if $c^{-1}x^sc=x^s$ then $(c^{-1}xc)^s=x^s$ so $c^{-1}xc=x$
and we conclude that $C_G(x^r)=C_G(x)$. In addition, as $G$ is
finitely generated and nilpotent it follows that torsion-free
Abelian subgroups of rank $1$ are infinite cyclic. Now, if $x,y\in
G$ such that $x$ is not a proper power and $x^n=y^m$, for some
$m,n\in\ZZ$, then $x$ and $y$ belong to a torsion-free Abelian
subgroup of rank $1$: namely the isolator of $x^n$, see \cite{K}.
Since this subgroup must be cyclic it follows that $y\in \la x
\ra$. Hence for all $y\in G$ there exists unique $x\in G$ with the
property that $x^n=y$ and whenever $y=z^m$ then $z\in\la x\ra$: we
call $x$ a  root of $y$ and say $x$ is a root element of $G$.
Since $c\ge 2$ and $G/\la e_c\ra$ is isomorphic to the semi-direct
product of the lattice $\ZZ^{c-1}$ and the infinite cyclic group
in the same way as $G$, the same properties hold in $G/\la e_c
\ra$. If $g\in G$ then we may choose $h\in G$ such that $h\la
e_c\ra$ is the root of $g\la e_c\ra$ in $G/\la e_c\ra$. Then
$g=h^nz$, for some $z\in \la e_c\ra$, so $C_G(g)=C_G(h^n)=C_G(h)$.
Thus we may assume that $t^r a$ is such that $t^r\la e_c \ra$ is a
root element of $G$. Clearly $C_G(t^r a)\supseteq \la t^ra,
e_c\ra$. Suppose that, for some $s\in \ZZ$ and $b\in A$, $t^sb\in
C_G(t^ra)$. We have $t^rat^sb=t^{r+s}(a\phi^s)b$ and
$t^sbt^ra=t^{r+s}(b\phi^r)a$ so it must be that
$(a\phi^s)b=(b\phi^r)a$. Writing $a=\sum_{i=1}^{c}\a_ie_i$ and
$b=\sum_{i=1}^{c}\b_ie_i$ (using additive notation for $A$) we
have
\begin{align*}
(a\phi^s)b& = (\a_1+\b_1)e_1\\
&+\left(\a_2+\b_2+{r\choose 1}\a_1\right)e_2 +\cdots\\[.5em]
&\cdots +\left(a_c+\b_c+{r\choose 1}\a_{c-1}+\cdots + {r\choose
c-1}\a_1\right)e_c,
\end{align*}
(where we take ${s\choose k}=0$ if $k>s$) and a similar expression
for $(b\phi^r)a$. Comparing coefficients of $e_i$'s in these two
expressions we see that for fixed $r$, $s$ and $a$ the elements
$\b_1,\ldots ,\b_{c-1}$ are uniquely determined. Hence, for each
$s\in\ZZ$ there is at most one coset $t^sb\la e_c\ra$ which is
contained in $C_G(t^ra)$. Now let $q=$lcm$(r,s)$, so there are
integers $u$ and $v$ such that $q=ur=vs$. Then $(t^ra)^u=t^qc\in
C_G(t^ra)$ and $(t^sb)^v=t^qd\in C_G(t^ra)$. From the above
$(t^ra)^u\la e_c\ra=(t^sb)^v\la e_c\ra$ and, since $t^ra\la
e_c\ra$ is a root element in $G/\la e_c\ra$, this means that
$t^sb\in \la t^ra, e_c\ra$. Therefore $C_G(t^ra)=\la t^ra,
e_c\ra$. The intersection of two such subgroups is $\la e_c\ra$
unless both subgroups are the same. It now follows that
$\cd(G)=2$, as claimed. 
\item\label{ex:nilp3} Let $G$ be the group
  constructed in the previous example, with $c \ge 3$, and let
  $H=G\ast_{Z(G)}G$.
  Then, by  Proposition \ref{prop:ofcd}, $\cd(G)=2$.
  Since $Z(H)=Z(G)$ we have $H/Z(H)\cong G/Z(G)\ast G/Z(G)$. Now
  $Z(G/Z(G))\cong \ZZ$ and so it follows from Proposition
  \ref{prop:ofcd}.\ref{i:cd3}
  that $\cd(H/Z(H)) =4$.
  This shows that, on taking the quotient of a group by its centre
  the centraliser dimension
  may increase. This example is directly comparable to 
  an example of R.~Bryant \cite{Bryant79} in which the original group $G$ has
  min-c but the factor group $G/Z(G)$ does not.
\end{enumerate}
\end{expl}

\affiliationone{%
   A.J. Duncan\\
   School~of~Mathematics~and~Statistics, University~of~Newcastle,
   Newcastle~upon~Tyne, NE1~7RU.\\
   \email{a.duncan@ncl.ac.uk}}
\affiliationtwo{% 
   I.V.~Kazatchkov and V.N.~Remeslennikov, \\
  Omsk~Branch~of~Mathematical~Institute~SB~RAS, 13~Pevtsova~Street,
  Omsk~644099,\\  Russia
   \email{kazatchkov@mail333.com\\
   remesl@iitam.omsk.net.ru}}
\end{document}